\documentclass[12pt,reqno,twoside]{amsart}

\usepackage{amsfonts,amsmath,amssymb}
\usepackage{mathrsfs,mathtools}
\usepackage{enumerate}
\usepackage{hyperref}
\usepackage{esint}
\usepackage{graphicx}
\usepackage{bm}
\usepackage{commath}
\usepackage{esint}
\DeclareMathAlphabet{\mathpzc}{OT1}{pzc}{m}{it}

\usepackage{placeins} 
\usepackage{flafter} 

\usepackage{caption}

\numberwithin{equation}{section}

\usepackage[dvips,bottom=1.4in,right=1in,top=1in, left=1in]{geometry}

\makeatletter
\def\eqnarray{\stepcounter{equation}\let\@currentlabel=\theequation
\global\@eqnswtrue
\tabskip\@centering\let\\=\@eqncr
$$\halign to \displaywidth\bgroup\hfil\global\@eqcnt\z@
  $\displaystyle\tabskip\z@{##}$&\global\@eqcnt\@ne
  \hfil$\displaystyle{{}##{}}$\hfil
  &\global\@eqcnt\tw@ $\displaystyle{##}$\hfil
  \tabskip\@centering&\llap{##}\tabskip\z@\cr}

\def\endeqnarray{\@@eqncr\egroup
      \global\advance\c@equation\m@ne$$\global\@ignoretrue}

\numberwithin{equation}{section}

\newcommand{\dt}[1]{\mathbf{DT}_{#1}}

\def\xvec{{\bf x}}

\def\dt{{\Delta t}}

\vfill

\title{Neural Network Representation of Time Integrators}

\thanks{This work is partially supported by NSF grant DMS-2110263 and the AirForce Office of Scientific Research under Award NO: FA9550-22-1-0248.}

\author{Rainald L\"ohner and Harbir Antil}
\address{R. L\"ohner and H. Antil. 
Center for Computational Fluid Dynamics and 
Center for Mathematics and Artificial Intelligence,
George Mason University,
Fairfax, VA 22030-4444, USA}


\subjclass[2010]{}

\begin{document}

\begin{abstract}
Deep neural network (DNN) architectures are constructed
that are the exact equivalent of explicit Runge-Kutta 
schemes for numerical time integration.
The network weights and biases are given,
i.e., no training is needed. In this way,
the only task left for physics-based integrators is
the DNN approximation of the right-hand side.
This allows to clearly delineate the approximation estimates
for right-hand side errors and time integration errors.
The architecture required for the integration of a simple 
mass-damper-stiffness case is included as an example.
\end{abstract}

\maketitle

\section{Introduction}

Considerable effort is currently being devoted to
neural nets, and in particular so-called deep neural nets (DNNs).
DNNs have been shown to be very good for sorting problems
(e.g. image recognition) or games (e.g. chess). Their use as ordinary
or partial differential equation (PDE) solvers is the subject of 
much speculation, with many variants such as Residual DNNs
\cite{EHaber_LRuthotto_2018a}, Physically Inspired NNs 
(PINNs) \cite{MRaissi_PPerdikaris_GEKarniadakis_2019a}, 
Numerically Inspired NNs (NINNs), Nudging NNs (NUNNs) 
\cite{HAntil_RLoehner_RPrice_2022a}, Fractional DNNs 
\cite{HAntil_RKhatri_RLohner_DVerma_2020a,HAntil_HCElman_AOnwunta_DVerma_2021a} 
and others being explored at present. 
The current situation is somewhat reminiscent of previous
attempts to use general, easy-to-use tools from other fields to
solve ordinary or partial differential equations. Examples include
the `discoveries' that one could use MS~Excel to solve ODEs 
\cite{Bilbao1996excel}, \cite{huseynov2013methodology}, 
Simulink for some simple PDEs 
\cite{RJGran_2007a,shampine1999solving}, cellular automata for ODEs and
PDEs \cite{wolfram2003new,wolfram1984cellular}, and ResNets for ODEs 
\cite{chen2018neural}. 

When solving time-dependent ODEs or PDEs, the resulting system is given
by:

$$ u_{,t} = r(u,t) ~~, $$

\noindent
where $u$ is the (scalar or vector) unknown, $r$ the (possibly nonlinear)
right-hand-side and $t$ time. The system can be integrated via explicit
Runge-Kutta (RK) schemes which are of the form:

\[
\begin{aligned}
 u^{n+1} &= u^{n} + \dt~b_i~r^i ~~,   \\
 r^i &= r(t^n+c_i \dt,~u^{n} + \dt~a_{ij}~r^j) 
                       ~~,~~ i=\overline{1,s},~~j=\overline{1,s-1} ~~ .  
\end{aligned}
\]                       
Any particular RK method is defined by the number of 
stages $s$ and the coefficients $a_{ij}, 1 \le j < i \le s$,
$b_i, i=\overline{1,s}$ and 
$c_i, i=\overline{2,s}$.  \\
Given that attempts are being made 
to replace time integrators by DNNs, one might ask: 
how should the architecture of the DNNs be in order to 
obtain the optimal time integration properties of RK schemes~?
This would clarify:
\begin{itemize}
\item[-] The minimum number of layers required for DNNs;
\item[-] The minimum width required for DNNs;
\item[-] The weights and biases required; 
\item[-] The overall efficiency of DNNs versus other alternatives; and
\item[-] Approximation properties of DNNs.
\end{itemize}

The remainder of the paper is organized as follows: Section~\ref{s:NN}
describes the standard neural network architectures. Section~\ref{s:poly} 
establishes that standard polynomials can be represented by the activation 
functions. Section~\ref{s:time} which illustrates how DNNs can be built
that result in standard time integrators.
The architecture required for the integration of a simple 
mass-damper-stiffness case is included as an example in Section~\ref{s:mass}.

\section{Neural Net Architectures} 
\label{s:NN}

A general DNN configuration 
consists of $L$ number of hidden layers along with one input and one
output layer. Each hidden layer, the input layer and the output layer
have $K$, $N$ and $J$ number of neurons, respectively.
The input-to-output sequence of such a DNN may be written as follows.

\begin{subequations}
\begin{align}
	\mbox{\em Input:} \quad
	G_{k}^{1} &= g \left( \sum_{n=1}^{N} w_{kn}^{1}M_{n}
                    +{[bias]}_{k}^{1} \right) , ~~
                                   k=\overline{1, K^{1}} \\
	\mbox{\em Hidden:} \quad                                   
	G_{k}^{l} &= g \left( \sum_{m=1}^{K^{l-1}} w_{km}^{l}G_{m}^{l-1}
                    +{[bias]}_{k}^{l} \right), ~~
                    k={}  \overline{1,K^{l}}, ~~  l={}  \overline{2,L}    \\        
	\mbox{\em Output:} \quad                    
	{BC}_{j} &=
			\phi \left( \sum_{m=1}^{K^L} w_{jm}^{L}G_{m}^{L} \right), 
    		   ~~ j= \overline{1,J},                                              
\end{align}
\end{subequations}
where $\phi(x), g(x)$ are activation functions, $w$ and $bias$ 
are the weights and biases.
Typical activation functions for $\phi(x), g(x)$ include:
\begin{itemize}
\item[-] Heaviside: $HS(x)=1$  for $x \ge 0, \phi(x)=0$ for $x < 0$
\item[-] Logistic: $LG(x)=1/(1+\exp(-x))$
\item[-] ReLU: $ReLU(x)=\max \{ 0,x \}$,
\item[-] HypTan: $HTAN(x)=\tanh(x)$.
\end{itemize}

Functions that do not have an `activation behaviour' but 
that have proven useful include:

\begin{itemize}
\item[-] Constant: $CO(x)=1$,
\item[-] Linear: $LI(x)=x$. 
\end{itemize}

\section{Polynomial Functions in 1-D}
\label{s:poly}
Let us now consider how to represent local polynomial 
functions via DNNs. An important question pertains to 
the activation functions used. In typical DNNs, these are
`switched on' when the input value crosses a threshold. 

\subsection{Constant Function}

Let us try to approximate the constant function $y(x)=1$ via
DNNs. The simplest way to accomplish this via `true activation
functions' with just one neuron would be via:
$$ DNN_c: \quad y(x)= HS(x-x_{\infty})  ~~, $$
where $x_{\infty}$ is a very large value. An alternative is to
use two neurons as follows:
$$ DNN_{c2}: \quad y(x)= HS(x+\epsilon) + HS(-x) ~~, $$
where $\epsilon$ would be of the order of machine roundoff.
Note that the desire to `activate' (which is seen as a
requirement of general DNNs) in this case has a negative
effect, prompting the need for either very large or small
numbers - something that may lead to slow convergence of
`learning' or numerical instabilities.
A far better alternative would have been the use of the
constant activation function $CO(x)$.

\subsection{Linear Function}

Let us try to approximate the linear function $y(x)=x$ via
DNNs and `true activation functions'. The obvious candidate
would be $ReLU(x)$. But as it has to work for all values of 
$x$ one can either use $ReLU(x)$ + $HS(x)$
$$ DNN_l:  \quad y(x)= ReLU(x-x_{\infty}) + x_{\infty} \cdot HS(x-x_{\infty})
                                                     ~~, $$
or:
$$ DNN_{l2}: \quad y(x)= ReLU(x+\epsilon) - ReLU(-x) ~~. $$
As before, the desire to `activate' has a negative effect,
prompting the need for either very large or small numbers.
A far better alternative would have been the use of the linear
activation function $LI(x)$.

In DNNs, one usually refrains from using higher order 
functions, trying to leverage the generality of lower order 
or differentiable activation functions.

\subsection{Higher Order Polynomial Functions in 1-D}

Consider now the polynomial 
$$ y(x) = a_j x^j ~~,~~ j \ge 2 ~~.  $$
The aim is to construct a DNN that would mirror this polynomial using 
the usual activation functions. Given that DNNs only act in
an additive manner, this is not possible. The usual recourse 
is to approximate it by a series of linear 
functions \cite{yarotsky2017error}. Another option is to transform to
logarithmic variables, add, and then transform back - but 
this would imply a major change in network architecture and
functions. Consider
$$ f(\xvec)= \sum_{j=1}^d a_j x_j ~~, $$
where $a_j$ are free coefficients and $x_j$ the spatial coordinates
in each dimension $j$. Note that only additions and `weights' ($a_j$)
are required, so the usual ReLU and HS functions should be able
to reproduce this function. But how many neurons are required~?
Borrowing from simplex (linear) finite element shape functions, 
one would have to build a linear function for each face of the ball
of elements (patch) surrounding a point. This implies a considerable number
of neurons for higher dimensional spaces. We refer to 
\cite{he2018relu,petersen2020neural,RDeVore_BHanin_GPetrova_2021a}.

\section{Explicit Timestepping for ODEs}
\label{s:time}

\begin{figure}[h!]
\includegraphics[width=6cm]{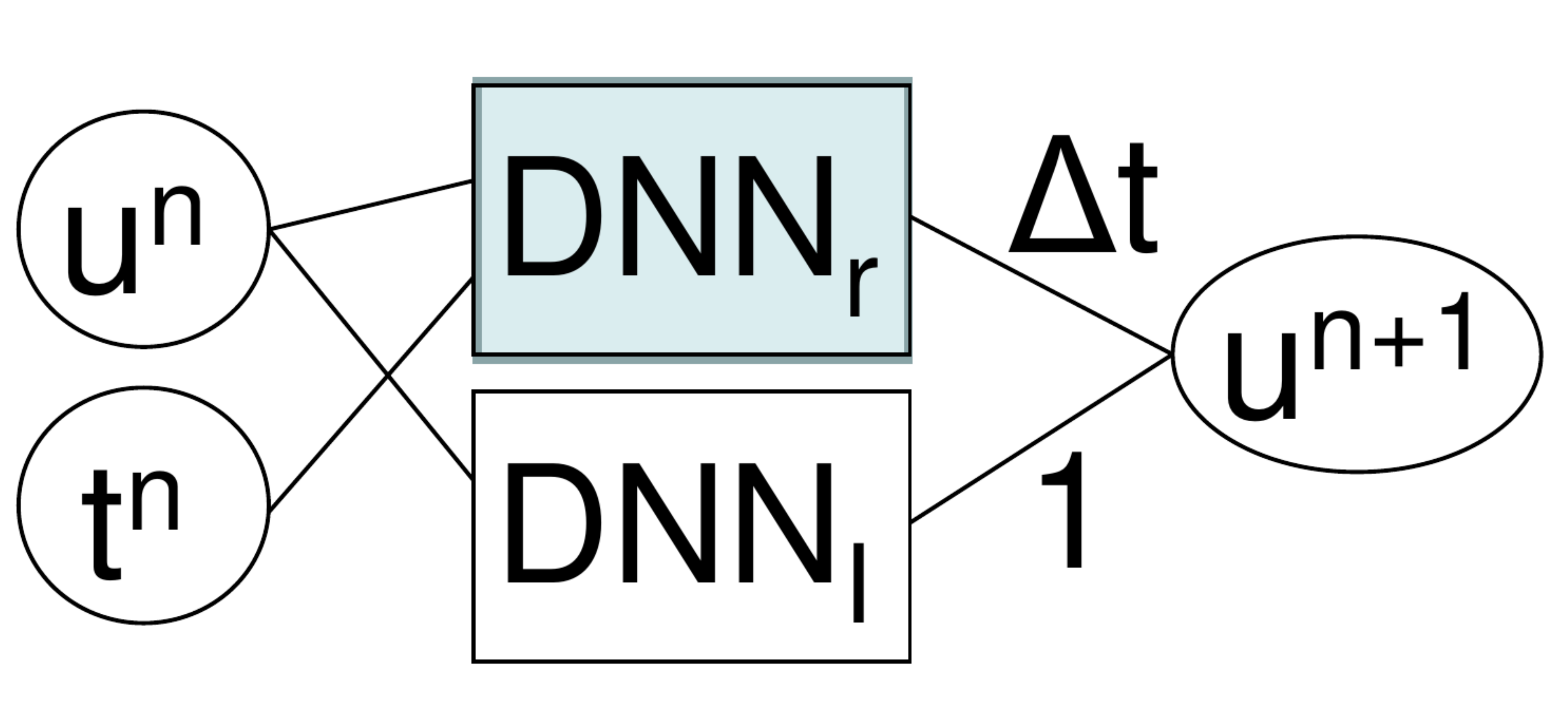}
\caption{\label{f:rk1}DNN for RK1 (Forward Euler) Scheme}
\end{figure}

Consider the typical scalar ODE of the form:
$$ u_{,t} = r(u,t) ~~. $$
Explicit time integration schemes take the right-hand-side $r$ at a known
time $t$ (or at several known times), and predict the unknown $u$ at some
time in the future based on it. The simplest such scheme is the 
forward Euler scheme, given by:
$$ u^{n+1} = u^{n} + \dt ~r(t^n, u^n) ~~. $$
Given that the function $r(u,t)$ is arbitrary, we will assume that
a DNN has been constructed for it. We will denote this approximation
of $r(u,t)$ as $DNN_r$. Note that in the scalar case this DNN has 
two inputs ($u,t$) and one output ($r(u,t)$). 
In order to obtain a complete DNN for the forward
Euler scheme, we need to enlarge $DNN_r$ by the `pass-through' value of
$u$. As was shown above, this can be accomplished with one layer
of 2~ReLU functions, or via one identity function.
We will denote this DNN as $DNN_I$ in the sequel.
The final DNN, shown in Figure~\ref{f:rk1} can then be denoted as:
$$ u^{n+1} = DNN_I(u^{n}) + \dt DNN_r(t^n, u^n) ~~. $$

In this and the subsequent figures we have highlighted the
`important' or 'essential' $DNN_r$ for $r(u,t)$. 
The generalization to higher order schemes is given by the
family of explicit {\bf Runge-Kutta} (RK) methods, which may be 
expressed as:
$$ 
\begin{aligned}
u^{n+1} &= u^{n} + \dt~b_i~r^i ~~,   \\
r^i &= r(t^n+c_i \dt,~u^{n} + \dt~a_{ij}~r^j) 
                       ~~,~~ i=\overline{1,s},~~j=\overline{1,s-1} ~~ .  
\end{aligned}                       
$$
Any particular RK method is defined by the number of stages $s$ and
the coefficients $a_{ij}, 1 \le j < i \le s$,
$b_i, i= \overline{1,s}$ and $c_i, i=\overline{2,s}$. 
These coefficients are usually arranged in a table known as a Butcher 
tableau (see Butcher~(2003)):
\begin{equation*}
  \begin{array}{c|ccccc}
    	         & r^1          & r^2  &   \dots    & r^{s-1} &         r^s         \\  
     \hline
    0           &                 & 		      &  	                      & \\
    c_2       & a_{21}      & 	              &  	    & 		     & \\
    c_3       & a_{31}      & a_{32}      &  	    & 		     & \\    
    \vdots    & \vdots      & \vdots     &  \ddots   & 		     & \\
    c_s       & a_{s1}      & a_{s2}   &   \dots    & a_{s,s-1} &  \\
    \hline
    \,        & b_1        & b_2       &  \dots     & b_{s-1} & b_s
  \end{array}
\end{equation*}
The two-step (2nd order) RK scheme is given by: 
\begin{equation*}
  \begin{array}{c|cc}
    	         & r^1          & r^2  \\
     \hline
    0           &                 & \\
    1/2        & 1/2           & \\
    \hline
    \,        & 0        & 1
  \end{array}
\end{equation*}

\begin{figure}[h!]
\includegraphics[width=8cm]{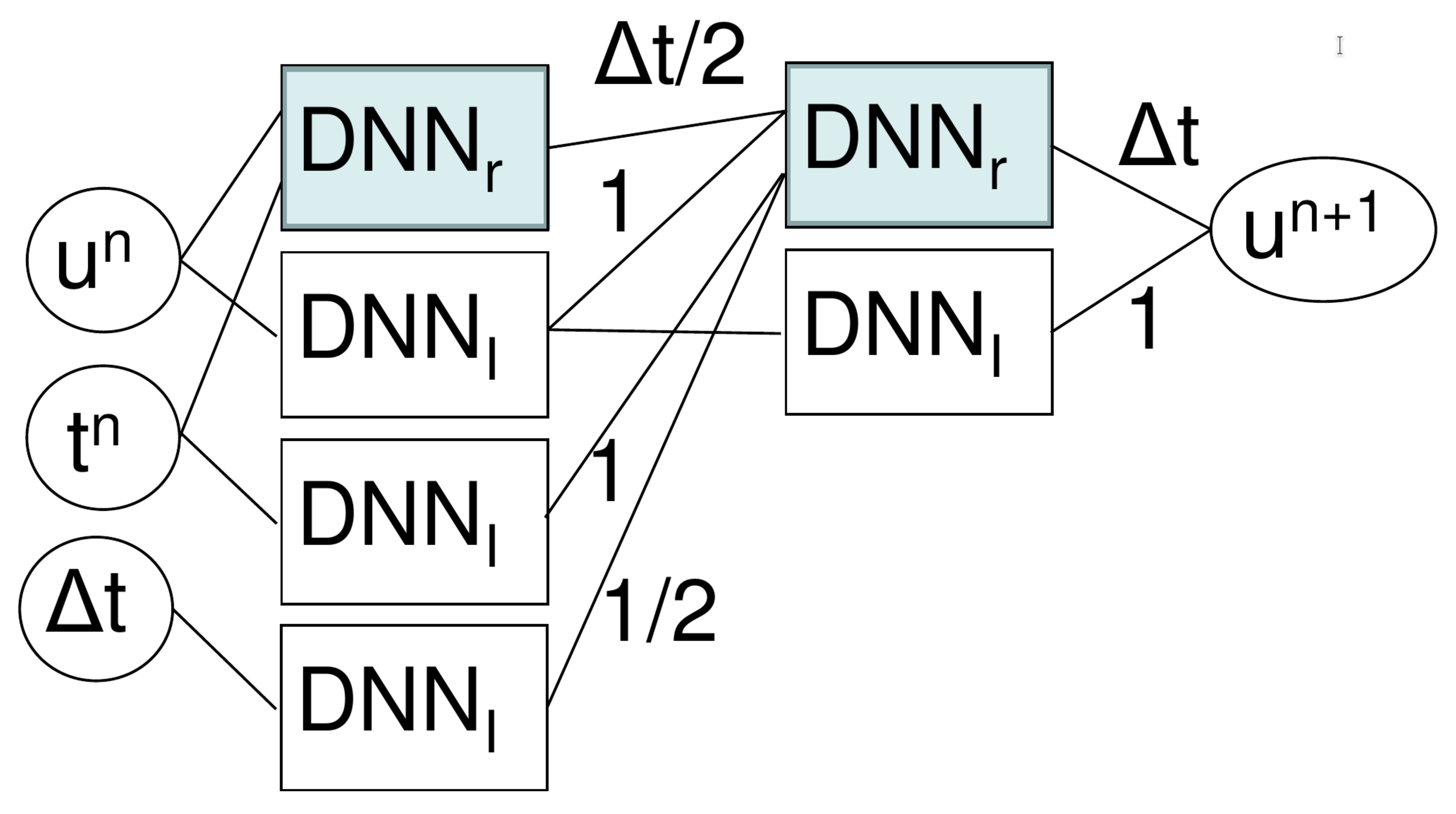}
\caption{\label{f:rk2}DNN for RK2 Scheme}
\end{figure}

For clarity, let us write the scheme out explicitly:
\begin{itemize}
\item[-] Step 1: $u^{n+1/2} = u^n + {\dt \over 2} r(u^n,t^n)$
\item[-] Step 2: $u^{n+1  } = u^n + \dt \, r(u^{n+1/2},t^{n+1/2})$. 
\end{itemize}
The DNN architecture required for this time integration scheme
is shown in Figure~\ref{f:rk2}.

The classic 4th order RK scheme is given by:
\begin{equation*}
  \begin{array}{c|cccc}
    	         & r^1          & r^2  &   r^3    & r^4 \\  
     \hline
    0          &                 & 		&  	                   \\
    1/2       & 1/2           & 	        &  	    & 		  \\
    1/2       & 0              & 1/2      &  	    & 		  \\    
    1          & 0              & 0         &   1     &  \\
    \hline
    \,        & 1/6             & 1/3       &  1/3     & 1/6
  \end{array}
\end{equation*}

\begin{figure}[h!]
\includegraphics[width=12cm]{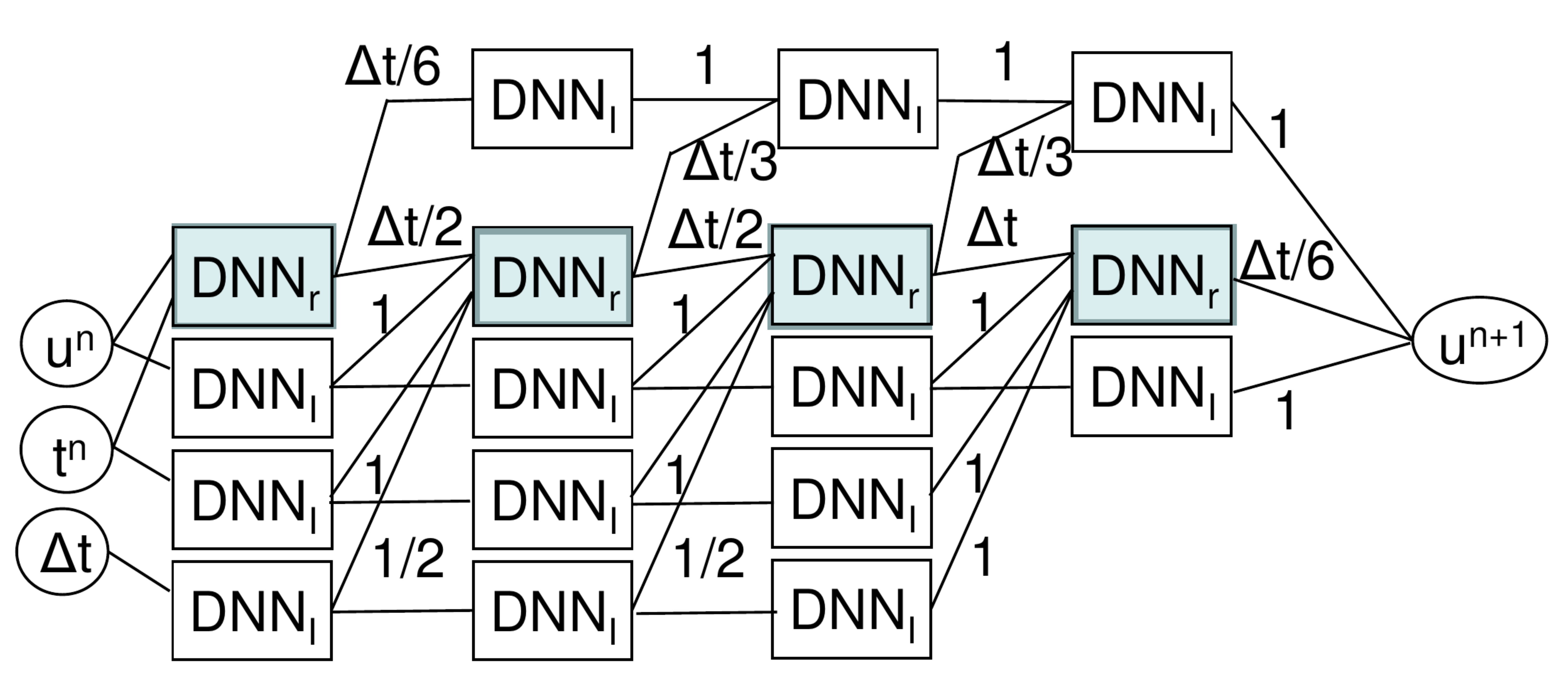}
\caption{\label{f:rk4}DNN for RK4 Scheme}
\end{figure}

For clarity, let us write the scheme out explicitly:
\begin{itemize}
\item[-] Step 1: $u^{n+1/4} = u^n + {\dt \over 2} r(u^n,t^n)$
\item[-] Step 2: $u^{n+1/3} = u^n + {\dt \over 2} r(u^{n+1/4},t^{n+1/2})$
\item[-] Step 3: $u^{n+1/2} = u^n +  \dt          r(u^{n+1/3},t^{n+1/2})$
\item[-] Step 4: $u^{n+1~~} = u^n + $
\item[ ] ${\dt \over 6}
\left[ r(u^n,t^n) + 2 r(u^{n+1/4},t^{n+1/2}) + 2 r(u^{n+1/3},t^{n+1/2})
     + r(u^{n+1/2},t^{n+1}) \right]$.
\end{itemize}
The DNN architecture required for this time integration scheme
is shown in Figure~\ref{f:rk4}.%

Observe that schemes of this kind require the storage of {\bf several}
copies of the unknown/right hand side, as the final result requires
$r^i, i=1,s$. Furthermore, as each right-hand side possibly requires the
information of all previous right-hand sides of the timestep, the
resulting neural net architecture deepens.

\section{Example: Mass-Damper-Stiffness System}
\label{s:mass}

Consider the simple mass-damper-stiffness system common to
structural mechanics, given by the scalar ODE:
$$ m \cdot x_{,tt} + d \cdot x_{,t} + c \cdot x = 0 ~~, $$
\noindent
where $m,d,c,x$ denote the mass, damping, stiffness and displacement
respectively. The ODE may be re-written as a first order ODE via:

$$ x_{,t} = v   ~~,~~ v_{,t} = - {d \over m} v - {c \over m} x ~~. $$

\begin{figure}[h!]
\includegraphics[width=12cm]{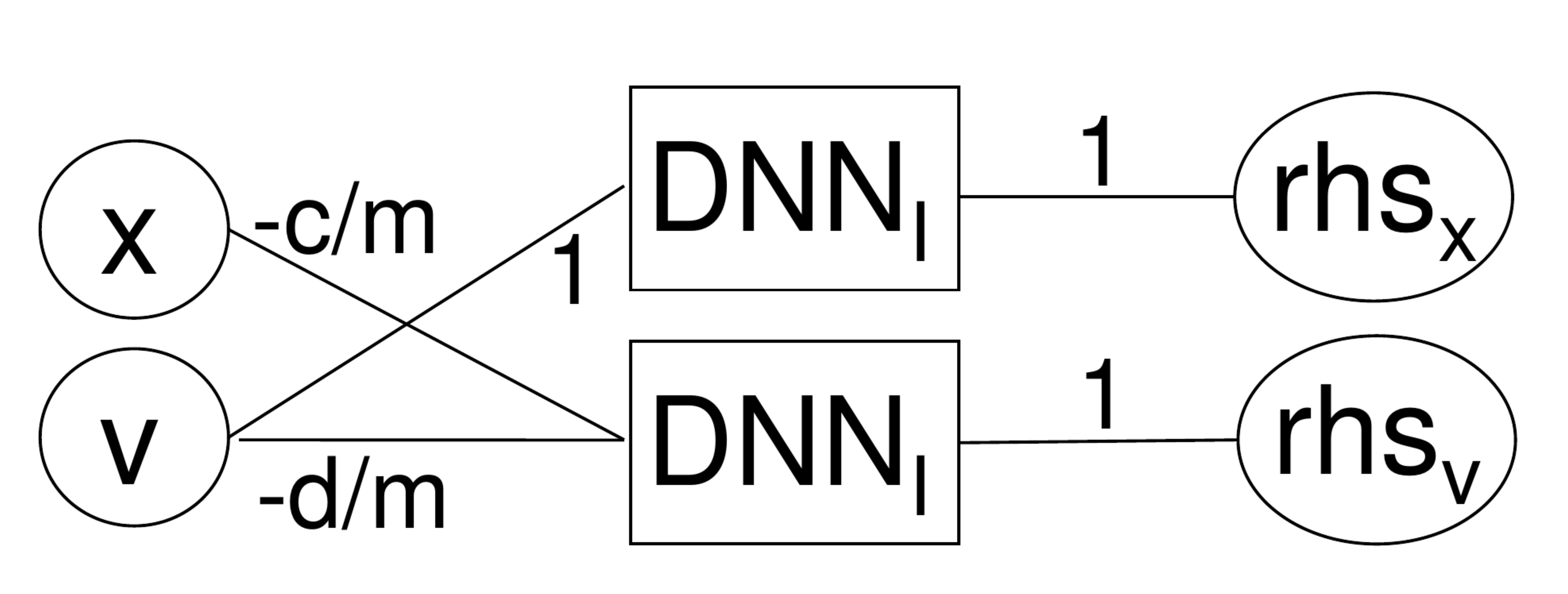}
\caption{\label{f:dnnMDC}DNN for Mass-Damper-Stiffness System}
\end{figure}

\noindent
The resulting $DNN_r$ is shown in Figure~\ref{f:dnnMDC}.

%
\section{Conclusions and Outlook}
Deep neural network (DNN) architectures were constructed
that are the exact equivalent of explicit Runge-Kutta
schemes for numerical time integration.
The network weights and biases are given,
i.e., no training is needed. In this way,
the only task left for physics-based integrators is
the DNN approximation of the right-hand side.
This allows to clearly delineate the approximation estimates
for right-hand side errors and time integration errors. \\
As the explicit Runge-Kutta schemes require the information of 
all previous right-hand sides of the timestep, the
resulting neural net architecture depth is proportional to
the number of stages - and hence to the integration order of 
the scheme. \\
As the DNN for the approximation of the right-hand side may
already be `deep', i.e. with several hidden layers, the final
DNN for high-order ODE integration many be considerable.

\bibliographystyle{plain}
\bibliography{refs}

\end{document}